\documentclass[11pt,reqno]{amsart}
\usepackage[T1]{fontenc}
\usepackage[margin=1in]{geometry}
\usepackage{amsmath,amssymb,amsthm}
\usepackage[hidelinks]{hyperref}
\hypersetup{pdftitle={An asymptotic resolution of the Gao--Thangadurai conjecture on the Harborth constant},pdfauthor={Guoqing Wang and Wenkai Yang}}
\newtheorem{theorem}{Theorem}[section]
\numberwithin{equation}{section}
\newtheorem*{conjecture}{Conjecture}
\newtheorem{lemma}[theorem]{Lemma}
\theoremstyle{remark}
\newtheorem*{remark}{Remark}
\newcommand{\Z}{\mathbb Z}
\newcommand{\stab}{\operatorname{Stab}}
\newcommand{\ex}{\operatorname{exp}}
\newcommand{\cF}{\mathcal F}
\newcommand{\gconst}{\mathsf g}
\title[The Gao--Thangadurai conjecture]{An asymptotic resolution of the Gao--Thangadurai\\
conjecture on the Harborth constant}

\author{Guoqing Wang \and Wenkai Yang}

\address{School of Mathematical Sciences \\ Tiangong University \\ Tianjin 300387, P.R. China}
\email{gqwang1979@aliyun.com}

\address{Center for Combinatorics \\ Nankai University, Tianjin 300071, P.R. China}
\email{9820250116@nankai.edu.cn}

\date{}
\begin{document}
\begin{abstract}
The Harborth constant $\gconst(G)$ of a finite abelian group $G$
of exponent $n$ is the least integer $k$ such that every subset
of $G$ of size at least $k$ contains $n$ elements whose sum is zero.
Let $C_n$ be a cyclic group of order $n$. In 2004, Gao and
Thangadurai \cite{gt} conjectured that $\gconst(C_n\oplus C_n)=2n+(-1)^n$ for all integers $n\ge3$.
We prove this conjecture for all sufficiently large integers $n$ and determine all
subsets of maximum size in $C_n^2$ that contain no $n$ elements with sum zero.
\end{abstract}
\maketitle

\section{Introduction}

The Erd\H os--Ginzburg--Ziv theorem \cite{egz} states that for every
positive integer $n$, any sequence of $2n-1$ integers contains $n$ terms
whose sum is divisible by $n$. Geometrically, these $n$ integers are
points on the line whose centroid is again an integer.

Kemnitz \cite{kemnitz83} asked the corresponding question in the plane.
How many points in $\Z^2$ suffice to guarantee that some $n$ of them
have a centroid in $\Z^2$? The centroid of $n$ points is a lattice point
precisely when both coordinate sums are divisible by $n$.
Reducing both coordinates modulo $n$ therefore gives an
equivalent additive problem in $C_n^2=C_n\oplus C_n$, where
$C_n=\Z/n\Z$ denotes the additive group of integers modulo $n$.
Each residue class has infinitely many lattice representatives, so
distinct lattice points may give repeated terms in $C_n^2$.

Taking $n-1$ copies of each vertex of the unit square shows that
$4n-4$ terms do not suffice. Kemnitz \cite{kemnitz83} conjectured
that the exact value is $4n-3$. Alon and Dubiner
\cite{ad93,ad} proved that $6n-5$ terms suffice for every $n$.
R\'onyai \cite{ronyai} reduced the upper bound to $4n-2$ for prime
$n$, and Gao \cite{gaoprimepower} extended this bound to prime powers.
Gao \cite{gaorestricted} also showed that, for each fixed $m\ge2$ and
all sufficiently large $n$, the conjecture for $n$ implies it for $mn$.
Thangadurai \cite{thangkemnitz} showed that $4n-3$ terms suffice
if some element occurs at least $\lfloor n/2\rfloor$ times.
Gao and Geroldinger \cite[Main Lemma~5.6]{ggcover}
showed that the conjecture would follow from a lower bound on the
number of affine hyperplanes avoiding the origin needed to cover certain
sets over prime fields.
Reiher \cite{reiher} proved the conjecture for every
positive integer $n$.

The corresponding problem for sets predates Kemnitz's
conjecture. In 1973 Harborth \cite{harborth} introduced the constant
$\gconst(C_n^d)$ for subsets of $C_n^d$. In the plane, how
many lattice points with pairwise distinct residue classes modulo $n$
suffice to guarantee that some $n$ of them have a lattice centroid?
In this case the sequence in $C_n^2$ is
\emph{squarefree}, so we may regard it as a set.
For a finite abelian group $G$ of exponent $n$, Gao and
Thangadurai \cite[Definition~1]{gt} define $\gconst(G)$ as the least
positive integer $k$ such that every subset of $G$ of size at least $k$
contains an $n$-element subset with sum zero. Following Marchan, Ordaz,
Ramos and Schmid \cite{mors}, we call $\gconst(G)$ the \emph{Harborth
constant} of $G$.

In $C_n$, the only subset of size $n$ is the whole group.
Its elements sum to $0$ if $n$ is odd and to $n/2$ if $n$ is even.
Hence $\gconst(C_n)=n$ for odd $n$ and $\gconst(C_n)=n+1$ for even $n$.
For $G=C_n^2$, Kemnitz
\cite{kemnitz83} proved that $\gconst(C_n^2)\ge 2n+(-1)^n$
for every $n\ge3$. Gao and Thangadurai \cite[Conjecture~1]{gt} conjectured that
these lower bounds are the precise values of $\gconst(C_n^2)$.

\begin{conjecture}
For every integer $n\ge3$,
\[
  \gconst(C_n\oplus C_n)=
  \begin{cases}
    2n-1, & n \text{ odd},\\
    2n+1, & n \text{ even}.
  \end{cases}
\]
\end{conjecture}

Harborth \cite{harborth} obtained $\gconst(C_3^2)=5$, and
Kemnitz \cite{kemnitz82} established $\gconst(C_n^2)=2n-1$ for $n=5,7$.
For an odd prime $p$, Sury and Thangadurai \cite[Proposition~4.6]{surythang}
proved that a set of $2p-1$ points in $C_p^2$ contains a $p$-point zero sum
if the multiplicities of its first coordinates are $(2,\ldots,2,1)$
or $(p-1,p-1,1)$. Gao and Thangadurai \cite{gt} proved $\gconst(C_p^2)=2p-1$
for every prime $p\ge67$ and also established $\gconst(C_4^2)=9$. Gao, Geroldinger
and Schmid \cite[Theorem~5.1]{ggs} extended the prime result to
$p\ge47$ and described the extremal sets of size $2p-2$.
Guillot, Marchan, Ordaz, Schmid and Zerdoum
\cite[Proposition~4.1]{guillot} established the next even case,
$\gconst(C_6^2)=13$. In a preprint, Graham \cite{graham} reports
$\gconst(C_p^2)=2p-1$ for every prime $p\ge11$.

Our first result covers all sufficiently large composite integers as well.

\begin{theorem}\label{main}
For every sufficiently large integer $n$,
$$
 \gconst(C_n\oplus C_n)=
 \begin{cases}
  2n-1,& n\ \text{odd},\\
  2n+1,& n\ \text{even}.
 \end{cases}
$$
\end{theorem}

For all sufficiently large $n$, we also determine the sets of size
$\gconst(C_n^2)-1$ that contain no $n$ elements with sum zero.
For odd $n$, the three forms in Theorem~\ref{inverse} agree with the
classification of Gao, Geroldinger and Schmid \cite[Theorem~5.1]{ggs}
for primes $p\ge47$. For even $n$, each extremal set is the union of
two cosets $a+H$ and $b+H$ of a cyclic direct summand $H$ of order $n$,
where $(b-a)+H$ generates $C_n^2/H$.

\begin{theorem}\label{inverse}
Let $n$ be sufficiently large and let $A\subseteq C_n^2$ have size
$2n+(-1)^n-1$. Then $A$ contains no $n$ elements with sum zero if
and only if, after a translation and a change of basis in $C_n^2$,
it has one of the following forms.

If $n$ is even, then $A=C_n\times\{0,1\}$.

If $n$ is odd, start with $C_n\times\{0,1\}$ and make one of the
following three changes. All horizontal coordinates belong to $C_n$.

\smallskip\noindent
{\rm (i)} Delete one point from each row.

\smallskip\noindent
{\rm (ii)} Delete $(x,0),(y,1),(y',1)$, where $y\ne y'$, and add
$(y+y'-x,2)$.

\smallskip\noindent
{\rm (iii)} Delete $(x,0),(x',0),(y,1),(y',1)$, where $x\ne x'$
and $y\ne y'$, and add $(u,-1)$ and $(v,2)$, where
\[
 u\in\{x+x'-y,x+x'-y'\},\qquad v\in\{y+y'-x,y+y'-x'\}.
\]
\end{theorem}

\section{Notation and preliminary results}

Throughout, $G$ is a finite abelian group with additive notation and
identity $0$. Let $C_n=\Z/n\Z$ and $C_n^2=C_n\oplus C_n$.
For an integer $d\ge1$, let $dG=\{dg:g\in G\}$. The exponent
$\ex(G)$ is the least positive $d$ such that $dG=\{0\}$.
For $H\le G$, let $\varphi_H:G\to G/H$ be the natural homomorphism and
let $[G:H]=|G/H|$ denote the index. If $X\subseteq G$, then
$X/H=\varphi_H(X)$.

We follow the sequence notation of Gao and Geroldinger
\cite[Section~2]{survey}. Let $\cF(G)$ denote the collection
of all finite unordered sequences of elements of $G$, with repetitions
allowed. For $S=g_1\cdots g_\ell\in\cF(G)$, its length is
$|S|=\ell$ and its sum is $\sigma(S)=g_1+\cdots+g_\ell$.
The notation $T\mid S$ means that $T$ is obtained by choosing some of the
positions of $S$. The sums of all subsequences of $S$ of length
$k$ form the set
$$
 \Sigma_k(S)=\{\sigma(T):T\mid S,\ |T|=k\}\qquad(0\le k\le |S|),
$$
with $\Sigma_0(S)=\{0\}$. A \emph{zero-sum subsequence of length
$k$} exists precisely when $0\in\Sigma_k(S)$. A sequence is
\emph{zero-sum free} if none of its nonempty subsequences has sum
zero. It is \emph{squarefree} if its terms are pairwise distinct.
We identify a finite set with the squarefree sequence of its elements.
When we project that sequence into $G/H$, we retain one term for
each original point. Thus two equal terms in the quotient can come
from distinct points of $G$.

For nonempty sets $X,Y\subseteq G$, let
$X+Y=\{x+y:x\in X,\ y\in Y\}$ and $X-Y=\{x-y:x\in X,\ y\in Y\}$.
In an expression such as $3X-2X$, the summands are chosen
independently from $X$ and may repeat. In $\Sigma_k(X)$, the $k$
summands are distinct. The \emph{stabilizer} of $X$ is the subgroup
$\stab X=\{g\in G:X+g=X\}$. We call $X$ \emph{aperiodic} if
$\stab X=\{0\}$. In $C_n$ every subgroup is $tC_n$ for a unique
divisor $t\mid n$. Its index is $t$, and its cosets are the
residue classes modulo $t$.

A \emph{primitive row} of $C_n^2$ is a coset of a cyclic direct
summand of order $n$. In coordinates adapted to that summand, the
parallel rows are the sets $C_n\times\{y\}$. We call $y$ the
\emph{height} and the first coordinate the \emph{horizontal coordinate}.
Translating all points by $g$ preserves every $n$-point sum, since
it changes that sum by $ng=0$.

A \emph{coset progression} with $r$ directions is a set of the form
$$
 x_0+H+\left\{\sum_{i=1}^r t_i g_i : t_i\in\Z,\ 0\le t_i<W_i\right\},
$$
where $H\le G$, $x_0,g_1,\ldots,g_r\in G$, and the $W_i$ are
positive integers. It is \emph{proper} if each of its elements has a
unique representation $x_0+h+\sum_i t_i g_i$ with $h\in H$ and
$0\le t_i<W_i$. When $H=\{0\}$, we call a proper progression a
\emph{box}. Its volume is the number of parameter choices,
$\prod_{i=1}^r W_i$.

We use Kneser's addition theorem in the form obtained by taking
$\ell=m$ in the theorem of DeVos, Goddyn and Mohar \cite[Theorem~1.3]{kneser}.
\begin{lemma}\label{kneser_bound}
Let $X_1,\ldots,X_k$ be nonempty finite subsets of $G$,
$X=X_1+\cdots+X_k$, and $H=\stab X$. Then
$|X/H|\ge 1+\sum_{i=1}^k (|X_i/H|-1)$.
\end{lemma}

For zero-sum subsequences of prescribed length, we use the
Erd\H os--Ginzburg--Ziv theorem \cite[p.~41]{egz} and
Reiher's solution of Kemnitz's conjecture \cite[Theorem~3.2 and the subsequent reduction]{reiher}.
\begin{lemma}\label{prescribed}
Let $b\ge1$. A sequence of at least $2b-1$ terms in $C_b$ contains a
zero-sum subsequence of length $b$. A sequence of at least $4b-3$
terms in $C_b^2$ contains a zero-sum subsequence of length $b$.
\end{lemma}

A quotient $G_1$ of $C_n^2$ with exponent $b$ embeds in $C_b^2$
by the invariant factor theorem.
Consequently, for every positive multiple $u$ of $b$,
\begin{equation}\label{quotient}
 |S|\ge u+3(b-1),\quad S\in\cF(G_1)\qquad\Longrightarrow\qquad
 0\in\Sigma_u(S).
\end{equation}
Indeed, successively remove $u/b$ disjoint zero-sum subsequences of
length $b$. Before the last choice at least $b+3(b-1)=4b-3$ terms
remain, so Lemma~\ref{prescribed} applies in $C_b^2$ at every step.
The chosen subsequences together have length $u$ and sum zero in $G_1$.

We use the following bound of Girard, Griffiths and Hamidoune
\cite[Theorem~1.1]{ggh}. Their exception requires a coset of an
elementary $2$-group of size at least four, which cannot occur in a
cyclic group.
\begin{lemma}\label{restricted}
Let $Y$ be a finite subset of a cyclic group. If $1\le k<|Y|$,
then $|\Sigma_k(Y)|\ge |Y|$.
\end{lemma}

The following two results of Savchev and Chen will be used for
sequences of heights and for their images in cyclic quotients.
\begin{lemma}\label{scfree}
Let $S$ be a zero-sum free sequence in $C_b$ with $|S|>b/2$.
There is a generator $g$ of $C_b$ such that, on expressing each term
as $ag$ with $1\le a<b$, the sum of the integers $a$ is less than $b$
\cite[Theorem~8]{sc}.
\end{lemma}
\begin{lemma}\label{scnlength}
Let $b/2<k<b$ and let $S$ be a sequence of length $b-1+k$ in $C_b$.
If $0\notin\Sigma_b(S)$, two distinct elements of $C_b$ have
combined multiplicity at least $2k$ in $S$ \cite[Corollary~7]{scn}.
\end{lemma}

To obtain a box from a set with many additive relations, we first
use the energy bound of Reiher and Schoen \cite[Theorem~1.2]{bsg}
with parameter $\varepsilon=1/4$. For a finite set $X\subseteq G$, let
$E(X)=|\{(x_1,x_2,x_3,x_4)\in X^4:x_1+x_2=x_3+x_4\}|$.
\begin{lemma}\label{energy_subset}
Let $K\ge1$ and let $X\subseteq G$ be finite and nonempty. If
$E(X)\ge |X|^3/K$, then there is a set $B\subseteq X$ such that
$|B|\ge 3|X|/(4\sqrt K)$ and $|B-B|\le 2^{51}K^4|B|$.
\end{lemma}

The Pl\"unnecke--Ruzsa inequalities and a result of Green and Ruzsa
allow us to pass from this difference set to a proper coset progression.
\begin{lemma}\label{pr_bound}
Let $X,Y$ be nonempty finite subsets of $G$ with $|X+Y|\le K|X|$.
For nonnegative integers $i,j$ with $i+j>1$,
$|iY-jY|\le K^{i+j}|X|$ \cite[Theorem~1.2]{petridis}, where
$0Y=\{0\}$.
\end{lemma}
\begin{lemma}\label{gr_box}
Let $B$ be a nonempty finite subset of $G$ with $|B+B|\le K|B|$.
Then $2B-2B$ contains a proper coset progression $P$ of rank at
most $2^9K^3\log_2(K+2)$ and size at least
$\exp(-2^{14}K^3\log_2^2(K+2))|B|$ \cite[Section~5, (5.1)]{gr}.
\end{lemma}

\section{A dense primitive row}

Throughout this section let $n\ge100$ and let $A\subseteq G=C_n^2$
satisfy $2n-2\le |A|\le2n+1$ and $0\notin\Sigma_n(A)$. We prove
that, for sufficiently large $n$, some primitive row contains more
than $4n/5$ points of $A$.

Lemma~\ref{layers} bounds the number of tuples when, for each sum
and each point, only a bounded number of tuples have that sum and contain that point.
\begin{lemma}\label{layers}
Let $X\subseteq A$, $|X|\le n$. Consider a family of ordered
$h$-tuples of distinct elements of $X$, where $h\ge2$. Suppose
that at most $J$ tuples with a given sum contain a given element
of $X$, with $J>0$. The family has fewer than $32hn^2J$ members.
\end{lemma}
\begin{proof}
We show that failure of the bound produces an $n$-term zero-sum
subsequence of $A$. Let $r(z)$ count the tuples with sum $z$, choose
$z_*$ with maximal $r(z_*)$, and set $\ell=\lfloor r(z_*)/(8hJ)\rfloor$. For $1\le j\le \ell$
let $Y_j=\{z:r(z)\ge8hJj\}$. If the asserted bound fails, then
$$
 \sum_j|Y_j|=\sum_{z\in G}\lfloor r(z)/(8hJ)\rfloor\ge3n^2.
$$
Also $hr(z_*)\le |X|J$ by incidence counting, so
$\ell\le n/(8h^2)$. For $U=\sum_jY_j$, let $H=\stab U$,
$w=[G:H]$, and $e=\ex(G/H)$. Lemma~\ref{kneser_bound} applied to $Y_1,\ldots,Y_\ell$ gives
$$
 w\ge 1-\ell+\sum_j |(Y_j+H)/H|\ge1-\ell+3w.
$$
Hence $e\le w\le(\ell-1)/2$. The set $U$ contains $\ell z_*+H$.
Round $\ell$ up to a multiple $\ell'$ of $e$, so $\ell'<3\ell/2$.
Since $e\mid\ell'$, we have $\ell'z_*\in H$. Adjoining
$\ell'-\ell$ copies of the singleton $\{z_*\}$ therefore produces
a sumset containing $H$.

Every element of $H$ is a sum of $h\ell'$ distinct points of $X$.
For each prescribed value in $Y_j$, choose a tuple in increasing order
of $j$. The earlier tuples use at most $h(j-1)$ points and thus
exclude fewer than $hjJ$ tuples at step $j$. For each additional
value $z_*$, fewer than $h\ell'J<8h\ell J$ tuples are excluded,
whereas $r(z_*)\ge8h\ell J$.
Since $h\ell'<3n/(16h)<n$ and $e$ divides both $n$ and $\ell'$,
the integer $u=n-h\ell'$ is positive and divisible by $e$. Since $h\ell'\ge2\ell$ and $e\le(\ell-1)/2$, we have
$h\ell'-3(e-1)\ge(\ell+9)/2>2$. Thus
$|A\setminus X|\ge n-2\ge u+3(e-1)$. By \eqref{quotient}, choose
$u$ points of $A\setminus X$ with sum in $H$, and cancel
that sum with the $h\ell'$ points inside $X$.
\end{proof}

\begin{lemma}\label{energy}
For $X\subseteq A$ of size $24\le m\le n$, one has
$E(X)>m^5/(4608n^2)$. Also, if $B_1,B_2\subseteq A$ are disjoint
of sizes $b,c\ge2$, with $b+c\le n$, then
$$
 \max_z |B_2\cap(z-B_1)|>\frac{b^2c^2}{1024n^2\max(b,c)}.
$$
\end{lemma}
\begin{proof}
We bound $E(X)$ by counting ordered triples with distinct entries.
Let $r_2(s)$ count the ordered pairs in $X^2$ with sum $s$, and
let $t=m^3/(576n^2)$. If at most $m^2/8$ pairs have
$r_2(s)>t$, then at least $m^3-3m^2-3m^3/8\ge m^3/2$ triples
have distinct entries and all three pair sums have multiplicity
at most $t$. Fixing a triple sum and one entry leaves at most
$3t$ such triples. Lemma~\ref{layers} with $h=3$ and $J=3t$
bounds their number strictly below $288n^2t=m^3/2$.
Thus $\sum_{s:r_2(s)>t}r_2(s)>m^2/8$, and
$E(X)=\sum_s r_2(s)^2>tm^2/8$.

To bound the intersection, consider tuples $(b_1,b_2,c_1,c_2)$
with $b_1,b_2\in B_1$ and $c_1,c_2\in B_2$, and unequal entries
within each pair. Since $B_1\cap B_2=\varnothing$, all four entries
are distinct. There are at least $b^2c^2/4$ such tuples.
For a fixed sum and a fixed element of $B_1$, there are two positions
for this element and $c$ choices for $c_1$.
The last $B_1+B_2$ pair has at most $\max_z|B_2\cap(z-B_1)|$
choices. Interchanging $B_1$ and $B_2$ gives the bound when the
fixed element belongs to $B_2$. Apply Lemma~\ref{layers} with
$X=B_1\cup B_2$, $h=4$ and $J=2\max(b,c)\max_z|B_2\cap(z-B_1)|$.
The condition $|X|\le n$ follows from $b+c\le n$.
\end{proof}

\begin{lemma}\label{struct}
Let $K=2^{240}$, $s=K^4$ and
$b_0=2^{-K^6}$. Then $A$ contains
$m$ points in a box
\[
 x_0+\Big\{\sum_{i=1}^r t_i g_i: 0\le t_i<W_i\Big\}
\]
with unique representations, where $r\le s$, $m\le n/2$,
$m\ge b_0 n$ and $m\ge V/K^5$, with $V=\prod_iW_i$.
\end{lemma}
\begin{proof}
We first find a subset with small doubling. Take $X\subseteq A$ of size $\lfloor n/2\rfloor$.
Since $|X|\ge n/3$, Lemma~\ref{energy} gives
$E(X)\ge |X|^3/K_0$, where $K_0=41472$.
Lemma~\ref{energy_subset} with parameter $K_0$ gives a set
$B\subseteq X$ satisfying
$$
 |B|\ge\frac{3|X|}{4\sqrt{K_0}}\ge\frac{n}{1024},\qquad
 |B-B|\le 2^{51}K_0^4|B|.
$$
Lemma~\ref{pr_bound} with $(X,Y)=(-B,B)$ and $(i,j)=(2,0)$ gives
$|B+B|\le (2^{51}K_0^4)^2|B|\le K|B|$.
Taking $(X,Y)=(B,B)$ and $(i,j)=(3,2)$ in Lemma~\ref{pr_bound}
now gives $|3B-2B|\le K^5|B|$.

To obtain the box, apply Lemma~\ref{gr_box} to $B$. There is a proper coset
progression $P\subseteq 2B-2B$ of rank at most
$2^9 K^3\log_2(K+2)$ and size at least
$\exp(-2^{14}K^3\log_2^2(K+2))|B|$. Averaging the identity
$$
 \sum_{z\in B-P}|B\cap(z+P)|=|B||P|
$$
over $B-P\subseteq 3B-2B$ gives a translate containing at least
$K^{-5}|P|$ points of $B$. This number is at least $b_0 n$,
since
$$
 2^{14}K^3\log_2^2(K+2)+\log(1024K^5)< K^6\log 2.
$$
The subgroup part of $P$ has at most two invariant factors.
Use one box coordinate for each such cyclic factor and translate the
origin of each parameter interval to zero. Properness gives unique
representations, so $V=|P|$, each $W_i\le n$, and
$r\le 2+2^9K^3\log_2(K+2)\le s$. Taking all the points of $B$ in
this translate gives $m\ge b_0n$ and $m\ge V/K^5$.
\end{proof}

To show that a dense primitive row must exist, we examine sets $A$ for
which every primitive row has at most $4n/5$ points. We first select
$n$ points whose horizontal sum vanishes modulo a divisor of $n$.

\begin{lemma}\label{project}
Suppose every primitive row contains at most $4n/5$ points of $A$.
Let $d\mid n$ and $n\ge200d^5$. After deleting at most $n/100$ points of $A$, one can
select $n$ remaining points with height sum zero in $C_n$ and first
coordinate sum zero modulo $d$.
\end{lemma}
\begin{proof}
We first prove the assertion for $d=1$. Let $S$ be the sequence of
remaining heights. If $|S|\ge2n-1$, Lemma~\ref{prescribed} gives
$0\in\Sigma_n(S)$ in $C_n$. Otherwise, for $k=|S|-n+1$, we have
$n-1-n/100\le k<n$ and $|S|=n-1+k$. Since $k>4n/5$,
Lemma~\ref{scnlength} with $b=n$ would give two heights of combined
multiplicity at least $2k>8n/5$ if the selection failed. Each of the
two rows contains at most $4n/5$ points.

For $d\ge2$, let $m=n/d$ and let $\ell\le n/100$ be the number
of points removed. A block consists of $d$
points whose two coordinate sums vanish modulo $d$. Its value in
$C_m$ is its height sum divided by $d$. It suffices to find $m$ disjoint
blocks whose values sum to zero in $C_m$, since their union has
$n$ points. We distinguish whether
there exist $k_0=\lceil(\ell+1)/d\rceil+5$ disjoint triples of
blocks with three different values in each triple.

If $k_0$ such triples exist, extend their blocks to a maximal
collection of disjoint blocks. At most $4d-4$ points remain:
otherwise Lemma~\ref{prescribed} in $C_d^2$ would give another block.
The number $b$ of blocks therefore satisfies
$b\ge(2n-2-\ell-4d+4)/d\ge2m-k_0$. If $b\ge2m-1$, Lemma~\ref{prescribed} in $C_m$ gives $m$ blocks
whose values sum to zero. Otherwise $m/2<b-m+1<m$, since $k_0\le m/100+7$ and
$m=n/d\ge200d^4$.
If no $m$ block values sum to zero, Lemma~\ref{scnlength} in
$C_m$ with $k=b-m+1$ gives two values accounting for all but at
most $2m-2-b\le k_0-2$ blocks. Each of the $k_0$ disjoint triples
contains at least one block with neither value, a contradiction.

If fewer than $k_0$ triples can be chosen, all but a bounded number
of remaining heights lie in at most two residue classes modulo $m$.
Discard a maximal collection of disjoint triples. From the remaining
points discard a maximal collection of disjoint blocks with pairwise
distinct values. This second collection has at most two blocks,
since three would form another triple. Let $U\subseteq C_m$ consist
of their values. In the two steps fewer than $3dk_0\le3\ell+3+18d$
points have been discarded. By maximality, any block formed from the
remaining points has value in $U$. All but at most $8d^2$ of their
heights belong to $U$ modulo $m$. Indeed, on the circle of circumference
$n$ mark every lift of $U$. Each of at most $2d$ open arcs
contains fewer than $4d-3$ points counted with multiplicity.
Otherwise Lemma~\ref{prescribed} in $C_d^2$ gives a block in that arc. Its mean
height, using integer representatives in an interval lifting the arc,
is an integer strictly between the endpoints. This integer modulo
$m$ is the block value and does not belong to $U$, a contradiction.
If $U$ is empty there is no block, so Lemma~\ref{prescribed} in
$C_d^2$ bounds the number of remaining points by $4d-4\le8d^2$.

To obtain enough copies of each label, keep only points whose
heights reduce to $U$ modulo $m$. For each residue represented
among these heights, choose an integer representative $a_i$.
Label each point $(x,y)$ in family $i$ by $z=(x\bmod d,z_2)\in G_1=C_d^2$,
where $y=a_i+mz_2$ in $C_n$. In each family discard the points
of each label occurring fewer than $d^3$ times. There are at most
$2d^2$ such groups of points, so fewer than $2d^5$ points are
discarded at this step. At least $2n-4\ell-5d^5$ points remain,
since $5+18d+8d^2+2d^5\le5d^5$ for $d\ge2$.
If one family has at least $n+d^3$ points, partition the points of
each label into $d$-tuples and a remainder of size at most $d-1$.
There are at most $d^2$ labels, so at most $d^2(d-1)<d^3$ points
belong to the remainders. Hence there are at least $n/d=m$ disjoint
$d$-tuples. The union of any $m$ of them has horizontal sum zero
modulo $d$ and height sum zero in $C_n$, since there are
$n$ terms $a_i+mz_2$ and each $z_2$ occurs a multiple of $d$ times.
Otherwise there are two families, each of size greater than $0.9n$.
Indeed $2n-4\ell-5d^5-(n+d^3)>0.9n$.

For each of the two families, indexed by $i=0,1$, we determine
all label sums of $k$ points for $n/3\le k\le2n/3$.
Fix a label $h_i$ in family $i$, and let $H_i\le G_1$ be generated
by the differences $h-h_i$ over all labels $h$ in that family.
For each $h\ne h_i$, reserve $\operatorname{ord}(h-h_i)-1$ disjoint
pairs of points with labels $h_i,h$. The total number $q$ of pairs
is less than $d^3$, and every label occurs at least $d^3$ times.
Selecting one point from each pair gives exactly $qh_i+H_i$, since
the pairs for each $h$ allow every multiple of $h-h_i$. Since
$q<d^3\le n/800$, we have $q\le k$ and $k+q<0.9n$.
There are thus at least $k-q$ unused points. Fixing $k-q$ of them
gives a coset of $H_i$ among the $k$-term sums. Every label belongs
to $h_i+H_i$, so this coset and the full set of $k$-term sums
are both $kh_i+H_i$.

The order of a single element in the quotient determines how many
points to select from each family. Let $H=H_0+H_1$ and
$g=(0,a_0-a_1)+m(h_0-h_1)\in G_1$. Let $t$ be the order of $g+H$
in $G_1/H$, so $t\mid d$. If $t<d$, set
$u=t\lfloor d/(2t)\rfloor$. Then $d/3\le u\le d/2$.
Select $um$ points from family $0$ and $(d-u)m$ from family $1$.
Both numbers lie in $[n/3,2n/3]$. Since $dm=n$ and $dG_1=\{0\}$,
the possible total label sums form the coset $um(h_0-h_1)+H$.
Since $ug\in H$, this coset contains $(0,-u(a_0-a_1))$.
Select points with this label sum. Their horizontal sum vanishes
modulo $d$, and their height sum in $C_n$ is
$um a_0+(d-u)m a_1-mu(a_0-a_1)=na_1=0$.

If $t=d$, each family lies on a primitive row. There is a
homomorphism from $G_1/H$ to $C_d$ taking $g+H$ to a generator.
Indeed, for every prime power $p^j$ exactly dividing $d$, some cyclic
factor of the $p$-primary part of $G_1/H$ contains a component of
$g+H$ of order $p^j$. That factor has order $p^j$ since $dG_1=\{0\}$.
The Chinese remainder theorem combines the projections onto these
factors into the required homomorphism. Composing with the quotient
map gives $\psi:G_1\to C_d$, $\psi(z_1,z_2)=cz_1+ez_2$, with
$H\subseteq\ker\psi$ and $\psi(g)$ a generator of $C_d$.
Let $e'\equiv e\pmod d$ and $e'\equiv1$ modulo each prime dividing
$m$ but not $d$. We show that $f:C_n^2\to C_n$,
$f(x,y)=cmx+e'y$, is surjective. If a prime $p$ divides both $m$
and $d$, then $\psi(g)\equiv e(a_0-a_1)\not\equiv0\pmod p$, so
$p\nmid e'$. If $p\mid d$ but $p\nmid m$, surjectivity of $\psi$
implies that at least one of $cm,e'$ is nonzero modulo $p$.
If $p\mid m$ but $p\nmid d$, then $e'\equiv1\pmod p$.
Thus $\gcd(cm,e',n)=1$, and there are integers $a,b$ with
$acm+be'\equiv1\pmod n$. The vectors $(e',-cm)$ and $(a,b)$ form a
basis of $C_n^2$ since their determinant is $1$ modulo $n$. Their
values under $f$ are $0$ and $1$, respectively, so $\ker f$ is the
cyclic direct summand generated by $(e',-cm)$. For each point $(x,y)$
of family $i$ with label $z$, we have
$f(x,y)=e'a_i+m\psi(z)=e'a_i+m\psi(h_i)$ since $z-h_i\in H_i\le H$.
Thus a coset of $\ker f$ contains the entire family. This primitive
row contains more than $0.9n>4n/5$ points of $A$, a contradiction.
\end{proof}

\begin{lemma}\label{rows}
Let $r_0\ge1$ be an integer and $c_0=2^{36}r_0^3$. If $n\ge200(c_0!)^5$
and every primitive row contains at most $4n/5$ points of $A$, then
more than $n/100$ points of $A$ lie in the complement of any union of
$r_0$ parallel primitive rows.
\end{lemma}
\begin{proof}
Suppose, to the contrary, that at most $n/100$ points of $A$ lie in
the complement of a union of $r_0$ parallel primitive rows. Choose
coordinates in which these rows are horizontal. We first select $n$
points with height sum zero and control their horizontal sum modulo
a divisor of $n$. Let $m_0=65536r_0$, $W=16r_0m_0$ and
$c_0=m_0W$, and set
$$
 a=\gcd\bigl(n,\operatorname{lcm}(1,\ldots,c_0)\bigr).
$$
Since $a\le c_0!$, our hypothesis gives $n\ge200a^5$, as required
by Lemma~\ref{project}. In every designated row
of size $c_y$ delete $\lfloor c_y/256\rfloor$ points, less than
$n/100$ in total since $|A|\le2n+1$. Apply Lemma~\ref{project}
with $d=a$ and these deletions. Of the $n$ selected points, at
least $n-n/100$ lie on the $r_0$ rows. Some such row therefore
has $k\ge(n-n/100)/r_0\ge n/(2r_0)$ selected points.
If this row originally contained $c$ points of $A$, then
$c-k\ge\lfloor c/256\rfloor\ge c/512$. Fix the selections on
all other rows and let $w\in C_n$ be the negative of their horizontal
sum. The $k$ selected points on the chosen row have a set of
horizontal coordinates $F\subseteq C_n$ with $\sigma(F)\equiv w\pmod a$.
We may vary this choice among all $c$ original points on that row.

We next preserve the horizontal sum modulo $a$ while retaining
enough independent choices on this row. Represent its $c$ original
horizontal coordinates in $[0,n-1]$.
Partition this interval into intervals of $W$ consecutive integers
(the last may be shorter). Within each form disjoint packets of
$m_0$ points until fewer than $m_0$ remain. At most
$(n/W+1)(m_0-1)\le n/(8r_0)$ points are left over, since $n\ge16r_0m_0$.
Since $c\ge n/(2r_0)$, we can retain $k_0=\lfloor n/(4096r_0m_0)\rfloor$
packets. Among the $c$ horizontal coordinates, choose a $k$-element
set $E\subseteq C_n$ containing exactly half of each retained packet.
This is possible since $k_0m_0/2\le\min(k,c-k)$. Pair each
$e\in E\setminus F$ with a distinct $f\in F\setminus E$.
The sum of the differences $f-e$ is
$\sigma(F)-\sigma(E)\equiv w-\sigma(E)\pmod a$. Repeatedly discard
nonempty zero-sum subsequences of these differences modulo $a$
until fewer than $a$ remain, using the partial-sum principle in $C_a$.
For the remaining pairs replace $e$ by $f$. The modified set has
horizontal sum congruent to $w$ modulo $a$, and at most $2a$
retained packets meet a modified pair.

To attain the exact horizontal sum, fix everything except the choices of $m_0/2$ in each undisturbed
packet. There are at least $k_0-2a\ge n/(8192r_0m_0)$ such packets.
Indeed $n\ge200(c_0!)^5\ge8192r_0m_0(2a+1)$. Let $tC_n$ be the stabilizer of the resulting
horizontal sumset. If $t>c_0$, each packet gives at least
$1+m_0^2/4$ sums modulo $t$. Indeed, order its points increasingly and move the selected $m_0/2$ indices from the first
half to the last by individual unit steps, always keeping them
distinct. The integer sums strictly increase over $m_0^2/4$ steps,
and their diameter is less than $m_0W$. In $C_t$ the total
sumset is aperiodic and proper. Lemma~\ref{kneser_bound} applied to the
packet sumsets modulo $t$ gives
$$
 t-2\ \ge\ \frac{n}{8192r_0m_0}\frac{m_0^2}{4}\ \ge\ 2n,
$$
which contradicts $t\mid n$. Thus $t\le c_0$, so $t\mid a$.
The modified set of horizontal coordinates has sum $z\equiv w\pmod a$.
Since $t\mid a$, $w-z\in tC_n$, so the horizontal sumset also contains
$w$. Together with the fixed points on the other rows this gives
$n$ points of $A$ whose two coordinate sums vanish, a contradiction.
\end{proof}

\begin{lemma}\label{singlerow}
For $0<a\le1$, define
$$
 b=\min(a/4,1/400),\quad r_a=1+\lceil3072/b^3\rceil,
 \quad c_a=2^{36}r_a^3,\quad F(a)=200(c_a!)^5 .
$$
If $n\ge F(a)$ and every primitive row contains at most $4n/5$ points
of $A$, then every primitive row contains
fewer than $a n$ points.
\end{lemma}
\begin{proof}
If a row contained at least $an$ points of $A$, we could cover all
but fewer than $n/100$ points with $r_a$ parallel primitive rows.
Mark the row and choose a subset $B_1$ of $\lfloor an/2\rfloor$
of its points. While at least $n/100$ points lie on unmarked parallel
rows, take $\lfloor n/200\rfloor$ of them as $B_2$. Then $B_1,B_2$
are disjoint, $|B_1|+|B_2|\le n$, and both sizes are at least $bn$
because $n\ge F(a)\ge\max(400,4/a)$. Lemma~\ref{energy} gives a
translate of $-B_1$ containing more than $b^3n/1024$ points of $B_2$.
The translate lies on one previously unmarked parallel row. Mark that row.
Since $|A|<3n$, fewer than $n/100$ points lie on unmarked rows
after at most $1+\lceil3072/b^3\rceil=r_a$ rows have been marked.
This contradicts Lemma~\ref{rows} with $r_0=r_a$ and $c_0=c_a$.
\end{proof}

The bound in Lemma~\ref{singlerow} will control the quotient by
the subgroup generated by the long directions of the box. For
the constants $K,s,b_0$ in Lemma~\ref{struct}, define
\begin{equation}\label{constants}
 D=16sK^5,\qquad \ell_0=D!,\qquad
 J=512s^2\ell_0K^{10},\qquad b_{i+1}=b_i^2/J\quad(0\le i<s).
\end{equation}
The integer $\ell_0$ is divisible by every possible stabilizer index
up to $D$. The numbers $b_i$ account for the points lost by
fixing short coordinates of the box.

\begin{lemma}\label{slicing}
Suppose $n>2s/b_s$. Some slice of the box in Lemma~\ref{struct}
has $r\ge1$ remaining coordinates and contains $m$ of the chosen
points. If $V$ denotes the volume of the slice, then
$$
 b_s n\le m\le n/2,\qquad m\ge V/K^5,\qquad
 W_i\ge Jn/m\quad(1\le i\le r).
$$
\end{lemma}
\begin{proof}
We successively fix short coordinates while preserving the density
and a lower bound for the number of points. After $i$ coordinates have been fixed, suppose the current box
contains at least $b_i n$ points and has density at least $K^{-5}$.
If a side has length $W<J/b_i$, fix the value of that coordinate
which retains the most points. The volume is divided by $W$, and the
new slice contains at least a $1/W$ fraction of the points.
Its density is therefore still at least $K^{-5}$, and its population
exceeds $b_i^2n/J=b_{i+1}n$. Repeat until all remaining sides have
length at least $J/b_i$. Fixing every coordinate would leave one site
containing at least $b_s n>1$ points, which is impossible.
At the stopping step, $m\ge b_i n\ge b_s n$, and every remaining
side has length at least $J/b_i\ge Jn/m$.
\end{proof}

We next bound the exponent of the quotient by the subgroup generated
by the remaining directions. Translate the surviving box so that its base point is zero. We translate
all of $A$, since the sum of $n$ points is translation invariant.
Let $K_1=\langle g_1,\ldots,g_r\rangle$, $v=\ex(G/K_1)$, and set
\begin{equation}\label{rowthreshold}
 a_0=\frac{b_s^2}{48\ell_0}.
\end{equation}
If $n\ge F(a_0)$ and every primitive row contains at most $4n/5$
points of $A$, then Lemma~\ref{singlerow} gives
\begin{equation}\label{smallv}
 v<\frac{m}{48\ell_0}.
\end{equation}
Indeed, Smith normal form for the inverse image of $K_1$ in $\Z^2$
gives a basis in which $K_1=uC_n\oplus vC_n$ with $u\mid v\mid n$.
The $m$ points lie on at most $n/v$ parallel primitive rows.
One row contains at least $mv/n$ points. If $v\ge m/(48\ell_0)$,
this is at least $m^2/(48\ell_0n)\ge b_s^2n/(48\ell_0)=a_0 n$,
a contradiction.
In particular $r\ge2$, since for cyclic $K_1$ the invariant factors
of $K_1$ have orders $n/v\mid n/u$, forcing $v=n$.

\begin{lemma}\label{axis}
Assume that the slice in Lemma~\ref{slicing} has $r\ge2$ remaining
coordinates. Partition its $m$ points into disjoint sets
$A_1,\ldots,A_r$ whose sizes differ by at most one. In $A_i$,
consider the lines
in the box obtained by varying only coordinate $i$, and let
$W=W_i$. Keep the lines containing at least $W/(4rK^5)$
points of $A_i$. For each retained line $L$, let $X_L$ be the set
of $i$-coordinates of the points of $A_i$ on $L$, and let
$c_L=|X_L|$. For any prescribed integers $k_L\in[c_L/3,c_L/2]$,
the stabilizer of $\sum_L\Sigma_{k_L}(X_L)$ has index at most $D$
in $C_n$, where the sum ranges over the retained lines in $A_i$.
\end{lemma}
\begin{proof}
We first estimate the total number of points on the retained lines.
Since $m>2s$, we have $|A_i|\ge m/(2r)$. There are $V/W$ lines
in direction $i$. The lines $L$ with $|A_i\cap L|<W/(4rK^5)$ contain at most $(V/W)W/(4rK^5)=V/(4rK^5)\le m/(4r)$
points of $A_i$ in total. The retained lines thus contain
$m_i\ge m/(4r)$ points. Because there is another side of length at least
$Jn/m$, we also have
\begin{equation}\label{massperwidth}
 m_i/W\ge V/(4rK^5W)\ge Jn/(4rK^5m)\ge48rK^5 .
\end{equation}
Let $t$ be the index of the stabilizer of $\sum_L\Sigma_{k_L}(X_L)$
in $C_n$. If $t>1$, the image of this sumset in $C_n/tC_n$ is
aperiodic and hence has at most $t-1$ elements. Lemma~\ref{kneser_bound} in $C_n/tC_n$ gives
\begin{equation}\label{axisgrowth}
 \sum_L\bigl(|\varphi_{tC_n}(\Sigma_{k_L}(X_L))|-1\bigr)\le t-2.
\end{equation}

Suppose $D<t\le W$. On a retained line $L$, the $c_L$ distinct
parameters in $[0,W-1]$ represent
$b=|\varphi_{tC_n}(X_L)|\ge c_L/\lceil W/t\rceil\ge c_Lt/(2W)$
residues modulo $t$. We show that the $k_L$-term sums represent at
least $b$ residues. This holds for $b=1$. For $b\ge2$, reserve one
point in each residue class and choose $j$ with
$1\le j<b$ and $0\le k_L-j\le c_L-b$. Such a choice exists since
$0<k_L<c_L$. Fix $k_L-j$ unreserved points and vary $j$ of the
representatives. Lemma~\ref{restricted}, applied to the $b$ residues
in $C_n/tC_n$, gives at least $b$ sums. Summing over the retained
lines and using $m_i\ge V/(4rK^5)$ and $t>D\ge16rK^5$, we obtain
\[
 \sum_L\bigl(|\varphi_{tC_n}(\Sigma_{k_L}(X_L))|-1\bigr)
 \ge \frac{m_it}{2W}-\frac VW\ge\frac{m_it}{4W}\ge t,
\]
where $m_i/W\ge48rK^5$ by \eqref{massperwidth}. This contradicts
\eqref{axisgrowth}.

Suppose $t>D$ and $t>W$. Let $h=\lfloor(t-1)/W\rfloor\ge1$. On
each retained line $L$, fix $k_L-j_L$ points, where $j_L=\min(k_L,h)$.
Among the $j_L+c_L-k_L$ remaining points, the $j_L$-term sums give
at least $1+j_L(c_L-k_L)$ residues modulo $t$. Indeed, as in the
proof of Lemma~\ref{rows}, increase the $j_L$ selected indices,
one unit at a time, from the first $j_L$ indices to the last $j_L$.
The resulting sums strictly increase over $j_L(c_L-k_L)$ steps
and have diameter at most $j_L(W-1)<t$.

To compare this bound with \eqref{axisgrowth}, let
$u_0=\sum_L k_L(c_L-k_L)$ and $k_*=\max_L k_L$.
Since $k_L(c_L-k_L)\ge c_L^2/6$ and $c_L\ge W/(4rK^5)$, we have
$u_0\ge Wm_i/(24rK^5)$. The bounds $m_i\ge m/(4r)$ and $W\ge Jn/m$
give $u_0\ge n$, and \eqref{massperwidth} gives $u_0\ge2W^2\ge2k_*W$.
If $h\ge k_*$, then $\sum_L j_L(c_L-k_L)=u_0\ge n\ge t$. If
$h<k_*$, then $j_L\ge hk_L/k_*$ on each line, so
$\sum_L j_L(c_L-k_L)\ge hu_0/k_*\ge2hW\ge t-1$,
since $t-1<(h+1)W$ and $h\ge1$. Each bound contradicts
\eqref{axisgrowth}. Hence $t\le D$.
\end{proof}

\begin{lemma}\label{dense_direction}
Suppose $n>2s/b_s$ and $n\ge F(a_0)$, with $s,b_s,a_0$ as in
\eqref{constants} and \eqref{rowthreshold}. If $A\subseteq C_n^2$,
$2n-2\le |A|\le2n+1$ and $0\notin\Sigma_n(A)$, some primitive row
contains more than $4n/5$ points of $A$.
\end{lemma}
\begin{proof}
Assume to the contrary that every primitive row contains at most $4n/5$
points of $A$. Let $B$ be the set of $m$ points in the slice of Lemma~\ref{slicing}.
By \eqref{smallv}, $v<m/(48\ell_0)<n$, so $r\ge2$. Let
$A_1,\ldots,A_r$ be the disjoint sets of Lemma~\ref{axis}.

We first find integers $k_L$ for which the total number of points is
divisible by $v\ell_0$. For each retained $i$-line $L$ in $A_i$ allow
$k_L$ to range over the multiples of $\ell_0$ in $[c_L/3,c_L/2]$.
Since $c_L\ge W_i/(4rK^5)\ge Jn/(4rK^5m)\ge24\ell_0$, this interval
contains such a multiple, and the largest and smallest allowed values
differ by at least $c_L/6-2\ell_0\ge c_L/12$. Adding these arithmetic
progressions over all retained lines gives every multiple of $\ell_0$
between the extreme totals. By the proof of Lemma~\ref{axis}, at least
$m/(4r)$ points of each $A_i$ lie on its retained lines. Thus the
interval of totals has width at least $m/48$, lower endpoint at least
$m/12$, and upper endpoint at most $m/2$. Since $v\ell_0<m/48$ by
\eqref{smallv}, this interval contains a multiple $\ell$ of $v\ell_0$.
Fix such $k_L$ with total $\ell$. In particular, $\ell\ge m/12>3v\ell_0$.

For each $i$ let $Y_i=\sum_L\Sigma_{k_L}(X_L)$, summed over the
retained $i$-lines in $A_i$. We show that $t_iC_n\subseteq Y_i$,
where $\stab Y_i=t_iC_n$. By Lemma~\ref{axis}, $t_i\le D$, so
$t_i\mid\ell_0$ and $t_i\mid k_L$. On each such line $L$ the
sequence $\prod_{x\in X_L}\varphi_{t_iC_n}(x)$ in
$C_n/t_iC_n\cong C_{t_i}$ has length $c_L\ge k_L+t_i-1$, since
$k_L\le c_L/2$ and $c_L\ge24\ell_0$. By Lemma~\ref{prescribed}
it contains $k_L/t_i$ disjoint zero-sum subsequences of length $t_i$,
since before the last one at least $c_L-k_L+t_i\ge2t_i-1$ terms remain.
Hence $0\in\varphi_{t_iC_n}(\Sigma_{k_L}(X_L))$ for each $L$,
and $0\in\varphi_{t_iC_n}(Y_i)$. Since $Y_i+t_iC_n=Y_i$, we have
$t_iC_n\subseteq Y_i$.

We now show that every element of
$$
 H'=\sum_{i=1}^r t_i\langle g_i\rangle
$$
is the sum of $\ell$ points of $B$. On a retained line $L$ in
$A_j$ with $j\ne i$, coordinate $i$ is constant. Since
$\ell_0\mid k_L$, the total fixed contribution to coordinate $i$
from all $A_j$ with $j\ne i$ belongs to $\ell_0C_n\subseteq t_iC_n$.
Adding this contribution to $Y_i$ still gives every element of
$t_iC_n$. The sets $A_i$ are disjoint, so the choices realizing
these coordinate sums are independent. The homomorphism
$(x_1,\ldots,x_r)\mapsto\sum_{i=1}^r x_i g_i$ maps
$\prod_{i=1}^r t_iC_n$ onto $H'$. Each choice uses exactly $k_L$
points on each retained line, so $H'\subseteq\Sigma_\ell(B)$.

It remains to find $u=n-\ell$ points in $A\setminus B$ whose sum
belongs to $H'$. Since $vG\subseteq K_1$ and $\ell_0K_1\subseteq H'$,
the exponent $e=\ex(G/H')$ divides both $v\ell_0$ and $n$. Thus
$u$ is a positive multiple of $e$, since $v\ell_0\mid\ell$ and
$\ell\le m/2<n$. Let $S=\prod_{a\in A\setminus B}\varphi_{H'}(a)$
be the sequence in $G/H'$, with one term for each point of $A\setminus B$.
Since $\ell\ge3v\ell_0\ge3e$ and $m\le n/2$, we have
$$
 u+3(e-1)\le n-3\le 2n-2-m\le |A|-m=|S|.
$$
By \eqref{quotient}, $0\in\Sigma_u(S)$, so there is a set
$T\subseteq A\setminus B$ of size $u$ with $\sigma(T)\in H'$.
Since $H'\subseteq\Sigma_\ell(B)$, there is a set $U\subseteq B$
of size $\ell$ with $\sigma(U)=-\sigma(T)$. Then $T\cup U\subseteq A$
has $u+\ell=n$ points and sum zero, a contradiction.
\end{proof}

\begin{remark}
Lemma~\ref{dense_direction}, and hence Theorems~\ref{main} and~\ref{inverse},
hold for every $n\ge n_0$, where one may take
$$
 n_0=2^{\,2^{\,2^{\,2^{1000}}}}.
$$
Let $x=s+3000$.
By \eqref{constants}, $s=2^{960}$, $D=2^{2164}$, and
$\log_2 \ell_0\le D\log_2 D<2^{2176}$. The recurrence gives
$$
 \log_2(1/b_s)
 \le 2^s\bigl(\log_2(1/b_0)+\log_2 J\bigr)
 <2^{s+2178},
$$
since $\log_2(1/b_0)=K^6=2^{1440}$ and
$\log_2 J=9+2\log_2 s+\log_2 \ell_0+10\log_2 K$.
In particular, $\log_2(2s/b_s)<2^x$ and
$\log_2(1/a_0)=\log_2(48\ell_0)+2\log_2(1/b_s)<2^x$.
Since $a_0<1/100$, the definition in Lemma~\ref{singlerow} gives
$r_{a_0}\le2^{19}a_0^{-3}$. Thus
$$
 \log_2 r_{a_0}<2^{x+3},\qquad
 \log_2 c_{a_0}=36+3\log_2 r_{a_0}<2^{x+6}.
$$
The elementary bound $\log_2(c_{a_0}!)\le c_{a_0}\log_2 c_{a_0}$ now gives
$\log_2 F(a_0)<2^{\,2^{x+8}}$.
Since $x+8<2^{1000}$, every $n\ge n_0$ satisfies both size conditions in Lemma~\ref{dense_direction}.

\end{remark}

\section{The Harborth constant and the extremal sets}

For a set $A$ with $0\notin\Sigma_n(A)$, the number
of points on a primitive row can force $|A|=2n-2$, as follows.

\begin{lemma}\label{partialrows}
Let $n\ge100$ and let $A\subseteq C_n^2$ satisfy
$2n-2\le |A|\le2n+1$ and $0\notin\Sigma_n(A)$. Suppose the largest
number of points of $A$ on a row in a fixed primitive direction is
strictly between $4n/5$ and $n$. Then $|A|=2n-2$. After a translation
and a change of height coordinate, the rows of heights $0,1$ contain
$n-r,n-t$ points of $A$, respectively, and one of the following holds:
$(r,t)=(1,1)$ and there are no other points, $(r,t)=(1,2)$ and there
is one further point of height $2$, or $(r,t)=(2,2)$ and there is
one further point of each height $-1$ and $2$.
\end{lemma}
\begin{proof}
We first show that a row containing the largest number of points of $A$
in the given direction attains all horizontal sums at each length
from $2$ to two less than that number. Translate this row to height
$0$, let $X=\{x\in C_n:(x,0)\in A\}$ and set $c=|X|=n-r$,
so $1\le r<n/5$. Since $\Sigma_{c-k}(X)=\sigma(X)-\Sigma_k(X)$,
it suffices to prove $\Sigma_k(X)=C_n$ for $2\le k\le c/2$.
For even $k$ choose two distinct terms with the prescribed sum.
For odd $k$ first choose one term and then choose two distinct
remaining terms with the adjusted sum. Such a pair exists because,
after deleting at most one term, the remaining set $Y$ satisfies
$|Y\cap(a-Y)|\ge2|Y|-n>2$ for every $a\in C_n$, whereas
$2x=a$ has at most two solutions. The set $X$ contains at least
$(2c-n-2)/2$ disjoint pairs $\{x,-x\}$ with $x\ne-x$.
At most three pairs meet the chosen terms, leaving at least
$c/4$ pairs since $3c\ge2n+16$. These pairs supply the remaining
terms without changing the horizontal sum.

Let $S$ be the sequence of nonzero heights of the points of $A$.
We show that some generator of $C_n$ occurs at least $n-5$ times in
$S$. We have $|S|=|A|-(n-r)\ge n+r-2$. If $T\mid S$ had
$\sigma(T)=0$ and $r+2\le |T|\le n-2$, then $2\le n-|T|\le c-2$.
Since $\Sigma_{n-|T|}(X)=C_n$, we could choose $n-|T|$ points at
height $0$ to cancel the horizontal sum of the points giving $T$,
contrary to $0\notin\Sigma_n(A)$. Suppose first that $S$ has a
minimal zero-sum subsequence $U$ of length $n-1$ or $n$. Deleting
one term from $U$ leaves a zero-sum free sequence of length greater
than $n/2$. By Lemma~\ref{scfree}, the positive representatives
of the remaining terms with respect to some generator $g$ sum to at
most $n-1$. Restoring the deleted term, with its representative in
$[1,n-1]$ relative to $g$, gives a positive multiple of $n$ less
than $2n$, since $\sigma(U)=0$. The representatives of $U$ therefore
sum to $n$, and at least $n-2$ of them are $1$.

If $S$ has no minimal zero-sum subsequence of length $n-1$ or $n$,
we obtain $n-5$ copies of a generator from a zero-sum free
subsequence. Remove disjoint minimal zero-sum subsequences from $S$
until the remainder is zero-sum free. Each removed subsequence has
length at most $r+1$: we have $0\notin\Sigma_j(S)$ for
$r+2\le j\le n-2$, and any $n$ terms in $C_n$ have a nonempty
zero-sum subsequence by the partial-sum argument. The total length
removed cannot first exceed $r+1$, since at that step the removed
terms would form a zero-sum subsequence of length in
$[r+2,2r+2]\subseteq[r+2,n-2]$. At least $n-3$ terms remain.
By Lemma~\ref{scfree}, their positive representatives with respect
to some generator $g$ sum to at most $n-1$. At least
$2(n-3)-(n-1)=n-5$ of these representatives are $1$.
In either case change the height coordinate from $y$ to $g^{-1}y$,
viewing $g$ as a unit modulo $n$. The row of height $1$ contains
$n-t$ points for some $r\le t\le5$, since no parallel row contains
more than $n-r$ points of $A$.

We restrict the heights of the $|A|-2n+r+t$ points not on the two
rows of heights $0,1$. Let $Y=\{x\in C_n:(x,1)\in A\}$.
For $1\le i<n-r$ and $1\le j<n-t$, Lemma~\ref{restricted} gives
$|\Sigma_i(X)|+|\Sigma_j(Y)|\ge(n-r)+(n-t)>n$. Thus
\begin{equation}\label{two-row-sums}
 \Sigma_i(X)+\Sigma_j(Y)=C_n\qquad
 (1\le i<n-r,\ 1\le j<n-t),
\end{equation}
since $\Sigma_i(X)\cap(a-\Sigma_j(Y))\ne\varnothing$ for every
$a\in C_n$. Suppose a point has height with integer representative
$h$ in $[t+1,n-r]$. Take this point, $h-1$ points of height $0$
and $n-h$ points of height $1$. Since $1\le h-1<n-r$ and
$1\le n-h<n-t$, we can make the horizontal sum zero by
\eqref{two-row-sums}. The integer representatives have sum $h+(n-h)=n$, so the height sum
is zero in $C_n$, contradicting $0\notin\Sigma_n(A)$.
Every remaining height thus has the form $-d$ with $1\le d\le r-1$
or $1+e$ with $1\le e\le t-1$.

The numbers of points of the two kinds are bounded by the numbers of
missing points on the rows of heights $0,1$. List the remaining
heights as $-d_1,\ldots,-d_a$ and $1+e_1,\ldots,1+e_b$, once per point.
We claim that
\begin{equation}\label{row-weights}
 \sum_{i=1}^a(d_i+1)\le r,\qquad \sum_{j=1}^b(e_j+1)\le t.
\end{equation}
If $\sum_i(d_i+1)>r$, take all $a$ points of heights $-d_i$,
together with $\sum_i d_i$ points of height $1$ and
$n-\sum_i(d_i+1)$ points of height $0$. Since
$a\le |A|-2n+r+t\le r+t+1\le11$ and $d_i\le4$, we have
$1\le\sum_i d_i\le44<n-t$ and $r<\sum_i(d_i+1)\le55<n$.
The $n$ chosen heights sum to zero, and \eqref{two-row-sums} allows
the horizontal sum to be zero as well, contrary to $0\notin\Sigma_n(A)$.
Hence $\sum_i(d_i+1)\le r$.
The change of height $y\mapsto1-y$ interchanges the two rows,
preserves $n$-point zero sums and sends $1+e_j$ to $-e_j$.
The same argument with $t$ in place of $r$ gives $\sum_j(e_j+1)\le t$.
Since each of the $a+b=|A|-2n+r+t$ points contributes at least $2$
to the sums in \eqref{row-weights}, we have
$2(|A|-2n+r+t)\le r+t$, or $r+t\le2(2n-|A|)$.

It follows that $|A|=2n-2$ and only three pairs $(r,t)$ are possible.
Indeed, $r,t\ge1$ gives $|A|\le2n-1$. If $|A|=2n-1$, then $r=t=1$
and a single point remains, although neither $1\le d_i\le r-1$
nor $1\le e_j\le t-1$ is possible. Thus $|A|=2n-2$ and $r+t\le4$.
The pair $(1,3)$ would give two points of heights $1+e_j$ with
$\sum_j(e_j+1)\ge4>t$, contrary to \eqref{row-weights}. Since
$1\le r\le t$, only $(1,1),(1,2),(2,2)$ remain. The first pair
leaves no further point. The second leaves one point, necessarily
of height $2$. For $(2,2)$, the two bounds in \eqref{row-weights}
force one further point of each height $-1$ and $2$.
\end{proof}

\begin{proof}[Proof of Theorem~\ref{main}]
For Kemnitz's lower bound \cite{kemnitz83}, let $B=C_n\times\{0,1\}$
if $n$ is even and $B=(C_n\setminus\{0\})\times\{0,1\}$ if $n$ is
odd. Then $|B|=2n+(-1)^n-1$. An $n$-point subset meeting both rows
has $j$ points of height $1$ for some $1\le j<n$, so its height sum
is nonzero. A single row contains $n$ points only when $n$ is even.
In this case its horizontal sum is $n(n-1)/2\equiv n/2\pmod n$.
Thus $0\notin\Sigma_n(B)$ and $\gconst(C_n^2)\ge 2n+(-1)^n$.

For the upper bound, assume to the contrary that $|A|=2n+(-1)^n$
and $0\notin\Sigma_n(A)$ for some $A\subseteq C_n^2$.
For sufficiently large $n$, Lemma~\ref{dense_direction} gives a row
with more than $4n/5$ points of $A$. Lemma~\ref{partialrows} then
forces a full row in that direction, since $|A|>2n-2$.
Translate the full row to height $0$. If $n$ is odd, the full row
sums to zero because $\sum_{x\in C_n}x=0$.

If $n$ is even, all $n+1$ points not on the full row would have to
lie on one parallel row. Let $S=h_1\cdots h_{n+1}$ be the sequence
of their heights. If $S$ had a zero-sum subsequence of length
$1\le j<n$, then $\Sigma_{n-j}(C_n)=C_n$ by Lemma~\ref{restricted},
so $n-j$ points of the full row could cancel the horizontal sum of
the $j$ corresponding points. Any $n$ of the $n+1$ heights have a
nonempty zero-sum subsequence by the partial-sum argument, and this
subsequence must therefore have length $n$. Omitting $h_i$ or $h_j$
gives $\sigma(S)-h_i=\sigma(S)-h_j=0$. Hence
$h_1=\cdots=h_{n+1}$. This would place $n+1$ distinct points on a row
of size $n$, a contradiction. Together with the lower bound, this gives
$\gconst(C_n^2)=2n+(-1)^n$.
\end{proof}

\begin{proof}[Proof of Theorem~\ref{inverse}]
Suppose $0\notin\Sigma_n(A)$. Lemma~\ref{dense_direction} gives a
primitive row with more than $4n/5$ points of $A$. If $n$ is even,
then $|A|=2n$, so Lemma~\ref{partialrows} forces a full row in this
direction. Translate the full row to height $0$ and let $S$ be the
sequence of heights of the other $n$ points. Since
$\Sigma_k(C_n)=C_n$ for $1\le k<n$ by Lemma~\ref{restricted},
the points giving any proper nonempty zero-sum subsequence of $S$
could be completed to an $n$-point zero sum with points from the full row. Thus
$S$ has no such subsequence. The partial-sum argument gives
$\sigma(S)=0$. Let $h\in C_n$ and $V\mid S$ satisfy $S=hV$. Then $V$
is zero-sum free and $|V|=n-1$. By Lemma~\ref{scfree}, the positive
representatives of $V$ with respect to some generator $g$ sum to at most
$n-1$. Thus every representative equals $1$ and $V=g^{n-1}$.
Since $h=-\sigma(V)=-(n-1)g=g$, the row of height $g$ is full.
Replacing $y$ by $g^{-1}y$ gives $A=C_n\times\{0,1\}$.
Conversely, a mixed selection of $n$ points from these two rows has
nonzero height sum. An entire row has horizontal sum
$n(n-1)/2\equiv n/2\pmod n$, so neither row gives an $n$-point zero sum.

If $n$ is odd, a full row would itself sum to zero since
$\sum_{a\in C_n}a=0$. Apply Lemma~\ref{partialrows} in the
direction given by Lemma~\ref{dense_direction}. Let $n-r,n-t$ be
the numbers of points on the rows of heights $0,1$, respectively.
If $(r,t)=(1,1)$, a choice of $n$ points must use between $1$ and
$n-1$ points of height $1$, so its height sum is nonzero.

If $(r,t)=(1,2)$, let $x$ be the missing horizontal coordinate
at height $0$ and $y,y'$ the two missing at height $1$.
An $n$-point choice of height sum zero must use the point $(v,2)$,
all $n-2$ points of height $1$, and one point of height $0$.
Since $\sum_{a\in C_n}a=0$, the required horizontal coordinate
at height $0$ is $y+y'-v$. Thus no $n$-point zero sum exists
if and only if $v=y+y'-x$.

If $(r,t)=(2,2)$, let $x,x'$ and $y,y'$ be the missing horizontal
coordinates at heights $0$ and $1$, respectively. An $n$-point
choice of height sum zero cannot omit both $(u,-1)$ and $(v,2)$,
since neither row is full. Nor can it use both, since that would
require $n-1$ points of height $1$. If it uses $(u,-1)$ alone
among these two points, it must use all $n-2$ points of height $0$
and one point of height $1$. The required horizontal coordinate
at height $1$ is $x+x'-u$. If it uses $(v,2)$ alone, the required
horizontal coordinate at height $0$ is $y+y'-v$. Thus no
$n$-point zero sum exists if and only if
$u\in\{x+x'-y,x+x'-y'\}$ and $v\in\{y+y'-x,y+y'-x'\}$.
\end{proof}

\section*{Declaration of generative AI and AI-assisted technologies}
The authors used a large language model during the preparation of this manuscript
for assistance with mathematical exploration, checking arguments,
organization and language editing. The authors are responsible for
the mathematical statements, proofs, references and final text.


\begin{thebibliography}{99}
\small
\setlength{\itemsep}{1pt}
\bibitem{ad93}
N. Alon and M. Dubiner, \emph{Zero-sum sets of prescribed
size}, in \emph{Combinatorics, Paul Erd\H os is Eighty}, Vol.~1,
J\'anos Bolyai Math. Soc., Budapest, 1993, 33--50.
\bibitem{ad}
N. Alon and M. Dubiner, \emph{A lattice point problem and
additive number theory}, Combinatorica \textbf{15} (1995), 301--309.
\bibitem{kneser}
M. DeVos, L. Goddyn and B. Mohar, \emph{A generalization of Kneser's
addition theorem}, Adv. Math. \textbf{220} (2009), 1531--1548.
\bibitem{egz}
P. Erd\H os, A. Ginzburg and A. Ziv, \emph{Theorem in the additive number
theory}, Bull. Res. Council Israel Sect. F \textbf{10} (1961), 41--43.
\bibitem{gaorestricted}
W. Gao, \emph{On zero-sum subsequences of restricted size},
J. Number Theory \textbf{61} (1996), 97--102.
\bibitem{gaoprimepower}
W. Gao, \emph{Note on a zero-sum problem},
J. Combin. Theory Ser. A \textbf{95} (2001), 387--389.
\bibitem{ggcover}
W. Gao and A. Geroldinger, \emph{Zero-sum problems and
coverings by proper cosets}, European J. Combin. \textbf{24} (2003), 531--549.
\bibitem{survey}
W. Gao and A. Geroldinger, \emph{Zero-sum problems in finite abelian
groups: a survey}, Expo. Math. \textbf{24} (2006), 337--369.
\bibitem{ggs}
W. Gao, A. Geroldinger and W. A. Schmid, \emph{Inverse zero-sum problems},
Acta Arith. \textbf{128} (2007), 245--279.
\bibitem{gt}
W. D. Gao and R. Thangadurai, \emph{A variant of Kemnitz Conjecture},
J. Combin. Theory Ser. A \textbf{107} (2004), 69--86.
\bibitem{ggh}
B. Girard, S. Griffiths and Y. O. Hamidoune, \emph{$k$-sums in
abelian groups}, Combin. Probab. Comput. \textbf{21} (2012), 582--596.
\bibitem{graham}
P. Graham, \emph{The Harborth constant of $C_p\oplus C_p$ for every
prime $p\ge11$}, Zenodo preprint, 4 November 2026,
\href{https://doi.org/10.5281/zenodo.23968472}{doi:10.5281/zenodo.23968472},
\href{https://aixiv.online/abs/2611.00487}{aiXiv:2611.00487}.
\bibitem{gr}
B. Green and I. Z. Ruzsa, \emph{Freiman's theorem in an arbitrary
abelian group}, J. Lond. Math. Soc. \textbf{75} (2007), 163--175.
\bibitem{guillot}
P. Guillot, L. E. Marchan, O. Ordaz, W. A. Schmid and H. Zerdoum,
\emph{On the Harborth constant of $C_3\oplus C_{3p}$},
J. Th\'eor. Nombres Bordeaux \textbf{31} (2019), 613--633.
\bibitem{harborth}
H. Harborth, \emph{Ein Extremalproblem f\"ur Gitterpunkte},
J. Reine Angew. Math. \textbf{262/263} (1973), 356--360.
\bibitem{kemnitz82}
A. Kemnitz, \emph{Extremalprobleme f\"ur Gitterpunkte}, doctoral
dissertation, Technische Universit\"at Braunschweig, 1982.
\bibitem{kemnitz83}
A. Kemnitz, \emph{On a lattice point problem}, Ars Combin.
\textbf{16-B} (1983), 151--160.
\bibitem{mors}
L. E. Marchan, O. Ordaz, D. Ramos and W. A. Schmid,
\emph{Some exact values of the Harborth constant and its plus-minus
weighted analogue}, Arch. Math. \textbf{101} (2013), 501--512.
\bibitem{petridis}
G. Petridis, \emph{New proofs of Pl\"unnecke-type estimates for product
sets in groups}, Combinatorica \textbf{32} (2012), 721--733.
\bibitem{reiher}
C. Reiher, \emph{On Kemnitz' conjecture concerning lattice-points in
the plane}, Ramanujan J. \textbf{13} (2007), 333--337.
\bibitem{bsg}
C. Reiher and T. Schoen,
\emph{Note on the theorem of Balog, Szemer\'edi, and Gowers},
Combinatorica \textbf{44} (2024), 691--698.
\bibitem{ronyai}
L. R\'onyai, \emph{On a conjecture of Kemnitz},
Combinatorica \textbf{20} (2000), 569--573.
\bibitem{sc}
S. Savchev and F. Chen, \emph{Long zero-free sequences in finite cyclic
groups}, Discrete Math. \textbf{307} (2007), 2671--2679.
\bibitem{scn}
S. Savchev and F. Chen, \emph{Long $n$-zero-free sequences in finite
cyclic groups}, Discrete Math. \textbf{308} (2008), 1--8.
\bibitem{surythang}
B. Sury and R. Thangadurai, \emph{Gao's conjecture on
zero-sum sequences}, Proc. Indian Acad. Sci. Math. Sci. \textbf{112}
(2002), 399--414.
\bibitem{thangkemnitz}
R. Thangadurai, \emph{On a conjecture of Kemnitz},
C. R. Math. Rep. Acad. Sci. Canada \textbf{23} (2001), 39--45.
\end{thebibliography}
\end{document}